\documentclass[reqno]{amsart}

\usepackage[english]{babel}
\usepackage{amsmath,amssymb,amsthm}
\usepackage{hyperref}
\usepackage[centering]{geometry}
\usepackage{xcolor}
\usepackage{color}
\renewcommand{\epsilon}{\varepsilon}            

\newtheorem{theorem}{Theorem}[section]   
\newtheorem*{theorem*}{Theorem}          

\theoremstyle{definition}

\newtheorem{remark}{Remark}[section]
\newtheorem*{acknow}{Acknowledgments}

\numberwithin{equation}{section}

\title[Isoparametric Hypersurfaces in product spaces of space forms]
{Isoparametric Hypersurfaces in product spaces of space forms}
\thanks{}

\author{Dong Gao\and Hui Ma\and Zeke Yao}

\address{D.~Gao, School of Science, Beijing
	University of Civil Engineering and Architecture, Beijing 102616, P.R. China}
\email{gaodong@bucea.edu.cn}

\address{H.~Ma, Department of Mathematical Sciences, Tsinghua
	University, Beijing, 100084, P.R. China}
\email{ma-h@mail.tsinghua.edu.cn}

\address{Z.~Yao, School of Mathematical Sciences, South China Normal University, Guangzhou 510631,
P.R. China}
\email{yaozkleon@163.com}


\date{}

\begin{document}
	
\begin{abstract}
We give a complete classification of isoparametric hypersurfaces in a product space
$M^2_{\kappa_1}\times M^2_{\kappa_2}$ of $2$-dimensional space forms for $\kappa_i\in \{-1,0,1\}$
with $\kappa_1\neq \kappa_2$. In fact we prove that any isoparametic hypersurface in
such a space has constant product angle function,
which enables us to remove the condition of constant principal curvatures
from the classification obtained recently by J.B.M.~dos Santos and J.P.~dos Santos.
\end{abstract}
	
\maketitle

\section{Introduction}\label{sect:1}

A hypersurface of a Riemannian manifold is called isoparametric if its nearby parallel
hypersurfaces have constant mean curvature. Isoparametric hypersurfaces provide the most
beautiful examples of hypersurfaces in Riemannian manifolds, and their study has received
a lot of attention. Isoparametric hypersurfaces in Euclidean and hyperbolic space
were classified due to the work of C. Somigliana, B. Segre and \'{E}. Cartan.
The classification of isoparametric hypersurfaces
in unit spheres is more complicated and was accomplished until recently.
We refer to the excellent survey in \cite{Chi,Tho} and the references therein.
	
However, we are still in the early stages of classifying the isoparametric hypersurfaces
in general Riemannian manifolds. According to \cite{GTY},
isoparametric hypersurfaces in $\mathbb{C}P^{2n}$ are homogeneous.
Later, a complete classification of isoparametric hypersurfaces in complex projective spaces
was obtained with the exception  of the case of  $\mathbb{C}P^{15}$ in \cite{D}.
And a classification of the isoparametric hypersurfaces in complex hyperbolic spaces
was given in \cite{DDS}. On the other hand, one can construct uncountably many isoparametric
families of hypersurfaces in Damek--Ricci spaces, and characterize those of them that have constant
principal curvatures by means of the new concept of
generalized K\"{a}hler angle \cite{DD}. Recently, a classification of isoparametric hypersurfaces
in any homogeneous $3$-manifold with $4$-dimensional isometry group was obtained in \cite{D-M}.

When the ambient manifold is a product space $M^2_{\kappa_1}\times M^2_{\kappa_2}$ of
$2$-dimensional space forms, there is a natural product structure $P$ on
$M^2_{\kappa_1}\times M^2_{\kappa_2}$ defined by $P(v_1,v_2)=(v_1,-v_2)$ for any tangent vector field
$(v_1, v_2)$. Here, $M^2_{\kappa_i}$ denotes the $2$-dimensional space form with constant sectional
curvature $\kappa_i$, for $\kappa_i\in \{-1,0,1\}$. Then given any orientable hypersurface $\Sigma$
of $M^2_{\kappa_1}\times M^2_{\kappa_2}$ with a unit normal vector field $N$, we can introduce an
important function $C$ defined by $C:=\langle PN,N\rangle$, where $\langle \cdot,\cdot\rangle$
denotes the standard product metric on $M^2_{\kappa_1}\times M^2_{\kappa_2}$.
Hereafter, for the sake of brevity, we call $C$ the product angle function of $\Sigma$.
Recently, Urbano \cite{Ur} classified the isoparametric hypersurfaces
in $\mathbb{S}^2 \times \mathbb{S}^2$. Subsequently, the authors \cite{GMY}
classified the isoparametric hypersurfaces in $\mathbb{H}^2 \times \mathbb{H}^2$.
If $\Sigma$ is an isoparametric hypersurface in $\mathbb{S}^2 \times \mathbb{R}^2$
with constant principal curvatures, Julio-Batalla \cite{Julio} showed that $\Sigma$
has constant product angle function $C$, and also classified the isoparametric hypersurfaces
with constant principal curvatures in $\mathbb{S}^2 \times \mathbb{R}^2$.
Later, J.B.M.~dos Santos and J.P.~dos Santos \cite{dos} generalized the result of
Julio-Batalla \cite{Julio} and showed that if $\Sigma$ is an isoparametric hypersurface in
$M_{\kappa_1}^2\times M_{\kappa_2}^2$ for $\kappa_1\neq \kappa_2$,
then the principal curvatures of $\Sigma$ are constant if and
only if the product angle function $C$ is constant. Furthermore, they classified the
isoparametric hypersurfaces with constant principal curvatures in the product spaces
$M_{\kappa_1}^2\times M_{\kappa_2}^2$ for $\kappa_1\neq \kappa_2$.

Recall that any isoparametric hypersurface in $M_{\kappa}^2\times M_{\kappa}^2$ ($\kappa=\pm 1$)
has constant product angle function $C$ (see \cite{Ur} and \cite{GMY}).
In this paper, we apply the Jacobi field computation appeared in \cite{dos} and the technique
developed by Urbano in \cite{Ur} to study isoparametric hypersurfaces in
$M_{\kappa_1}^2\times M_{\kappa_2}^2$ for $\kappa_1\neq \kappa_2$.
Combining with the results in \cite{Ur} and \cite{GMY}, we obtain our first result.

\begin{theorem}\label{thm:1.1}
Any isoparametric hypersurface in $M_{\kappa_1}^2\times M_{\kappa_2}^2$
has constant product angle function $C$, where $\kappa_i\in \{-1,0,1\}$ and
$\kappa_1^2+\kappa_2^2\neq0$.
\end{theorem}

\begin{remark}\label{rek:1.1}
Instead of an auxiliary function in \cite{dos} by means of the mean curvature of
the nearby parallel hypersurfaces and the determinant of the matrix obtained from
the Jacobi vector field, we focus on this determinant function itself in our proof.
The fact that its derivatives are constant provides us enough equations to deduce
that the product angle function $C$ is constant.
\end{remark}

Then, taking account of Theorem 1 of \cite{dos}, we know that the isoparametric hypersurfaces
in $M_{\kappa_1}^2\times M_{\kappa_2}^2$ for $\kappa_1\neq \kappa_2$
have constant principal curvatures. Thus, as a direct application of Theorem 2 of \cite{dos},
we obtain the following classification  of isoparametric hypersurfaces in
$M_{\kappa_1}^2\times M_{\kappa_2}^2$ for $\kappa_1\neq \kappa_2$.

\begin{theorem}\label{thm:1.2}
Let $\Sigma$ be an isoparametric hypersurface in $M_{\kappa_1}^2\times M_{\kappa_2}^2$
for $\kappa_i\in \{-1,0,1\}$ and $\kappa_1\neq \kappa_2$. Then, up to isometries,
$\Sigma$ is an open part of one of the following hypersurfaces:
\begin{itemize}
\item[(1)]
$\Gamma_1 \times  M_{\kappa_2}^2$, where $\Gamma_1$ is a  curve with constant curvature
in $M_{\kappa_1}^2$;

\item[(2)] $M_{\kappa_1}^2 \times \Gamma_2$, where $\Gamma_2$ is a  curve with constant
curvature  in $M_{\kappa_2}^2$;

\item[(3)] $\Psi(\mathbb{R}^3) \subset \mathbb{H}^2 \times \mathbb{R}^2$,
where $\Psi:\mathbb{R}^3\rightarrow\mathbb{H}^2\times\mathbb{R}^2$,
$(t,r,s)\rightarrow(p(t,r),q(t,s))$ is an immersion given by
			\begin{equation*}
				\begin{aligned}
					&p(t,r)=\cosh(t\sqrt{c})\gamma_1(r)+\sinh(t\sqrt{c})N(r),\\
					&q(t,s)=\gamma_2(s)+t \sqrt{1-c} W_0,
				\end{aligned}	
			\end{equation*}
and $\gamma_1(r)=\{(\frac{2+r^2}{2},r,\frac{r^2}{2})\in \mathbb{H}^2|\ r\in (-\infty,\infty)\}$
is a horocycle, $N(r)=(\frac{r^2}{2},r,\frac{-2+r^2}{2})$
is a unit normal vector field to $\gamma_1(r)$, $\gamma_2(s)=X_0+ s V_0$, $V_0$ and $W_0$ are constant
orthogonal unit vectors in $\mathbb{R}^2$, $X_0$ is a vector in $\mathbb{R}^2$, $0<c<1$ is a constant.
\end{itemize}
\end{theorem}

\begin{remark}\label{rek:1.2}
The expression of $\Psi(\mathbb{R}^3)$ in (3) of Theorem \ref{thm:1.2}
is slightly different from the one in Theorem 2 of \cite{dos}.
In fact, $\Psi(\mathbb{R}^3)$ provides a hypersurface which is the union of a family of
parallel surfaces given by the products of horocycles in $\mathbb{H}^2$
and straight lines in $\mathbb{R}^2$.
\end{remark}
		
\section{Preliminaries}\label{sect:2}
		
We denote by $M_{\kappa}^2$, the Euclidean plane $\mathbb{R}^2$ endowed with the canonical
metric $\langle \cdot,\cdot \rangle$ of constant sectional curvature $0$ when $\kappa=0$,
and the unit $2$-sphere $\mathbb{S}^2=\{x\in \mathbb{R}^3\mid x_1^2+x_2^2+x_3^2=1\}$
endowed with the canonical metric $\langle \cdot,\cdot \rangle$ of constant
sectional curvature $1$ when $\kappa=1$, and the hyperbolic plane
$\mathbb{H}^2=\{(x_1,x_2,x_3)\in \mathbb{R}_1^3\mid -x_1^2+x_2^2+x_3^2=-1,x_1>0\}$
endowed with the canonical metric $\langle \cdot,\cdot \rangle$ of
constant sectional curvature $-1$ when $\kappa=-1$,
where $\mathbb{R}_1^3$ is the three-dimensional Minkowski
space with the Lorentzian metric.
		
The complex structure $J$ on $M_0^2$ is defined by
$$
J_p(x_1,x_2)=(-x_2,x_1),
$$
for any $p\in M_0^2$ and $(x_1,x_2)\in T_pM_0^2$. The complex structure $J$ on $M_{\kappa}^2$
for $\kappa=1$ or $-1$ is defined by
$$
J_pv=\begin{cases}
p \wedge v, & \text{if } \kappa=1,\\
p \boxtimes v, &  \text{if } \kappa=-1,
\end{cases}
$$
for any $p\in M_{\kappa}^2$ and $v\in T_pM_{\kappa}^2$, where
$\wedge$ stands for the cross product in $\mathbb{R}^3$ and
$\boxtimes$ denotes the Lorentzian cross product defined by
$$
(a_1,a_2,a_3) \boxtimes (b_1,b_2,b_3)=
(a_3b_2-a_2b_3,a_3b_1-a_1b_3,a_1b_2-a_2b_1).
$$
Then we have two natural complex structures on $M_{\kappa_1}^2\times M_{\kappa_2}^2$ given by
$$
J_1=(J, J), \ \ J_2=(J,-J).
$$
The product structure $P$ on $M_{\kappa_1}^2\times M_{\kappa_2}^2$ is defined by
$$
P(X_1,X_2)=(X_1,-X_2), \ \ \forall X_1\in TM_{\kappa_1}^2,X_2\in TM_{\kappa_2}^2.
$$
It is easy to see that $P=-J_1J_2=-J_2J_1$, $P^2=\text{Id}$ and
\begin{equation}\nonumber
\langle PX,Y \rangle=\langle PY,X \rangle, \ \
\forall X,Y\in T(M_{\kappa_1}^2\times M_{\kappa_2}^2).
\end{equation}
Moreover, $\bar{\nabla}P=0$, where $\bar{\nabla}$ is the
Levi-Civita connection on $M_{\kappa_1}^2\times M_{\kappa_2}^2$.
		
The curvature tensor $\bar{R}$ of $M_{\kappa_1}^2\times M_{\kappa_2}^2$
with respect to the Riemannian product metric is given by
\begin{equation}\label{eqn:2.1}
\begin{aligned}
\bar{R}(X, Y, Z, W)= & \frac{\kappa_1}{4}\{\langle X, P W+W\rangle\langle Y, P Z+Z\rangle-\langle X,
P Z+Z\rangle\langle Y, PW+W\rangle\} \\
& +\frac{\kappa_2}{4}\{\langle X, P W-W\rangle\langle Y, P Z-Z\rangle-\langle X, P Z-Z\rangle\langle Y,
PW-W\rangle\},
\end{aligned}
\end{equation}
where $X, Y, Z, W \in T(M_{\kappa_1}^2\times M_{\kappa_2}^2)$.
	

Let $\Sigma$ be an orientable real hypersurface of $M_{\kappa_1}^2\times M_{\kappa_2}^2$
with the unit normal vector field $N$. Then, with respect to the product structure $P$,
the  product angle function $C$ of $\Sigma$  and
a vector field $V$ tangent to $\Sigma$ are defined by
\begin{align*}
C&:=\langle PN,N \rangle=\langle J_1N,J_2N \rangle,\\
V&:=PN-CN.\nonumber
\end{align*}
It is clear that $-1\leq C \leq 1$ and $||V||^2:=\langle V,V\rangle=1-C^2$. From the Guass equation
and \eqref{eqn:2.1} we can get the Ricci curvature of $\Sigma$
\begin{equation*}
\begin{aligned}
\operatorname{Ric}(X,X)= & \frac{\kappa_1}{4}[(1-C)\langle X, X\rangle+(1-C)\langle X,
PX\rangle+\langle PX, N\rangle^2 ] \\
& +\frac{\kappa_2}{4}[(1+C)\langle X, X\rangle-(1+C)\langle X, PX\rangle+\langle PX, N\rangle^2 ]\\
&+H\langle AX, X\rangle-\langle AX, AX\rangle,
\end{aligned}
\end{equation*}
where $X\in T\Sigma $. So the scalar curvature of $\Sigma$ is given by
\begin{equation}\label{eqn:2.2}
\rho=\kappa_1(1-C)+\kappa_2(1+C)+H^2-||A||^2.
\end{equation}
		
\section{Proof of Theorems \ref{thm:1.1} and \ref{thm:1.2}}\label{sect:6}

In this section, we are mainly devoted to the proof of Theorem \ref{thm:1.1}.
Theorem \ref{thm:1.2} follows as a direct consequence of Theorem \ref{thm:1.1} and
Theorems $1$ and $2$ of \cite{dos}.

Suppose that $\Sigma$ is an isoparametric hypersurface of $M_{\kappa_1}^2\times M_{\kappa_2}^2$ for
$\kappa_1\neq \kappa_2$. Let $\Phi=(\phi,\psi): \Sigma\rightarrow M_{\kappa_1}^2\times M_{\kappa_2}^2$
denote the inclusion map, and $N=(N_1,N_2)$ be the unit normal vector field of $\Sigma$.
Consider an open subset $\mathcal{U}$ of $\Sigma$ defined by
$\mathcal{U}=\{p \in \Sigma \mid C^{2}(p)<1\}$.
If $\mathcal{U}$ is empty, then $C^{2}\equiv 1$ holds on $\Sigma$.
By \cite{dos}, we know that $\Sigma$ is an open part of $\Gamma_1 \times  M_{\kappa_2}^2$ or
$M_{\kappa_1}^2 \times \Gamma_2$, where $\Gamma_i$ is a curve with constant curvature
in $M_{\kappa_i}^2$, for $i=1$ or $2$.
		
Now we suppose that $\mathcal{U}$ is not empty. Then the nearby parallel
hypersurfaces to any given connected open subset of $\mathcal{U}$, still denoted by $\mathcal{U}$,
are given by $\Phi_{l}: \mathcal{U} \rightarrow M_{\kappa_1}^2\times M_{\kappa_2}^2$, where
$$
\Phi_l(p)=\left(\exp _{\phi(p)} (l N_{1}(\phi(p))),
\widetilde{\exp}_{\psi(p)} (l N_{2}(\psi(p)))\right),
\quad p\in\mathcal{U},
\quad l \in(-\epsilon, \epsilon), \quad \Phi_{0}=\Phi,
$$
$\exp$ and $\widetilde{\exp}$ denote the exponential maps of $M_{\kappa_1}^2$
and $M_{\kappa_2}^2$, respectively.
		
For $p\in \mathcal{U}$, let $\gamma_{p}(l)$ be the geodesic with $\gamma_p(0)=p$ and
$\gamma^{\prime}_p(0)=N_p$.
Then the unit normal vector field  of the parallel hypersurface $\Phi_l(\mathcal{U})$ is given by
$N_{\gamma_{p}(l)}^l=\gamma^{\prime}_p(l)$.
Using the fact that $P$ is parallel, we get
\begin{equation*}
C_l=\langle PN^l, N^l\rangle=C,
\quad \forall l \in(-\epsilon, \epsilon).
\end{equation*}
		
Next, we can choose on $\Phi_l(\mathcal{U})$ the following orthonormal frame
$$
E_{1}^l=\frac{V_l}{\sqrt{1-C^{2}}}, \quad E_{2}^l=\frac{J_{1} N^l+J_{2} N^l}{\sqrt{2(1+C)}},
\quad E_{3}^l=\frac{J_{1} N^l-J_{2} N^l}{\sqrt{2(1-C)}} ,
$$
where $V_l=PN^l-CN^l$. Thus it can be verified  that  $\{E_{i}^l\}_{i=1}^3$ are parallel
along  the geodesic $\gamma_{p}(l)$.
		
Let $E_i^0 = E_i$ and $A_{ij}=\langle AE_i, E_j\rangle$ the components of the shape operator $A$
with respect to the unit normal vector field $N$ on $\mathcal{U}$.
Then,
$(\Phi_l)_*{E_i}$ is a Jacobi field along $\gamma_{p}(l)$ with initial conditions
$$
((\Phi_l)_*{E_i})|_{l=0}= E_i,\quad ((\Phi_l)_*{E_i})^{\prime}|_{l=0}=-AE_i=-\sum_{j=1}^{3}A_{ji}E_j.
$$
By solving the Jacobi equations  (see sect. 10.2.2 of \cite{B-C-O}, or see \cite{dos}), we obtain that
\begin{equation*}
\begin{aligned}
(\Phi_l)_*{E_i}=&\left(\delta_{1 i}-l A_{1 i}\right) E_{1}^l+\left(\delta_{2 i}
C_{\delta_1}(l)-A_{2 i} S_{\delta_1}(l)\right) E_{2}^l
+\left(\delta_{3 i} C_{\delta_2}(l)-A_{3 i}
S_{\delta_2}(l)\right) E_{3}^l,
\end{aligned}
\end{equation*}
where $\delta_{ij}$ is the usual Kronecker symbol for $i, j=1, 2$ or $3$ and  the functions
$S_{\delta_a}(l)$ and $C_{\delta_a}(l)$ are given by
\begin{equation*}
S_{\delta_a}(l)=\begin{cases}
l &   \text { if } \delta_a=0, \\
\frac{1}{\sqrt{-\delta_a}} \sinh \left(l \sqrt{-\delta_a}\right) & \text { if } \delta_a<0, \\
\frac{1}{\sqrt{\delta_a}} \sin \left(l\sqrt{\delta_a}\right) & \text { if } \delta_a>0,
\end{cases}
\quad \,
C_{\delta_a}(l)= \begin{cases} 1 & \text { if } \delta_a=0, \\
\cosh \left(l \sqrt{-\delta_a}\right) & \text { if } \delta_a<0, \\
\cos \left(l \sqrt{\delta_a}\right) & \text { if } \delta_a>0,\end{cases}
\end{equation*}
for $a \in\{1,2\}$, with $\delta_1=\frac{\kappa_1(1+C)}{2}$ and $\delta_2=\frac{\kappa_2(1-C)}{2}$.
Notice that
\begin{equation}\label{eqn:3.1}
S_{\delta_a}^{\prime}(l)=C_{\delta_a}(l),  \quad \quad  C_{\delta_a}^{\prime}(l)=-\delta_a
S_{\delta_a}(l).
\end{equation}
		
Denote $(\Phi_l)_*{E_i}= \sum_{j=1}^3 Q_{ji} E_j^l$,
where the matrix $Q=(Q_{ji})$ is given by
\begin{equation*}
Q=
\left(
\begin{array}{ccc}
1- lA_{11} & - l A_{12}  &    - l A_{13}\\
-A_{12}S_{\delta_1}(l) & C_{\delta_1}(l)-A_{22}S_{\delta_1}(l) & -A_{23}S_{\delta_1}(l) \\
-A_{13}S_{\delta_2}(l)  & -A_{23}S_{\delta_2}(l) & C_{\delta_2}(l)-A_{33}S_{\delta_2}(l) \\
\end{array}
\right).
\end{equation*}
			
Furthermore, according to Theorem 10.2.1 of \cite{B-C-O}, we obtain
\begin{equation*}
A_l((\Phi_l)_*{E_i})=-\sum_{j=1}^3 Q_{ji}^{\prime} E_{j}^l,
\end{equation*}
where $A_l$ denotes the shape operator of  the hypersurface $\Phi_l(\mathcal{U})$
associated with the unit normal vector field $N^l$.				
Thus it follows that $A_l=-Q^{\prime} Q^{-1}$ and the mean curvature
of the nearby hypersurface $\Phi_l(\mathcal{U})$ is given by
\begin{equation*}
H(l)=-{\rm tr}\left(Q^{\prime} Q^{-1}\right)
=-\frac{({\rm det} Q)^{\prime}}{{\rm det} Q}.
\end{equation*}
It implies that $({\rm det} Q)^{\prime}=-H(l) {\rm det} Q$
holds on  $\mathcal{U}$.
From $Q(0)={\rm I d}$, an inductive argument shows that
$$
\left(\frac{d^{k} {\rm det} Q}{d l^{k}}\right)(0), \quad k \geq 0,
$$
are constant on $\mathcal{U}$.
By straightforward computations, we have
$$
\begin{aligned}
\operatorname{det} Q=&(1-lA_{11})C_{\delta_1}(l)C_{\delta_2}(l)
+(-A_{22}+lH_{12})S_{\delta_1}(l) C_{\delta_2}(l) \\
&+(-A_{33}+lH_{13})C_{\delta_1}(l) S_{\delta_2}(l)
+(H_{23}-lK)S_{\delta_1}(l) S_{\delta_2}(l),
\end{aligned}
$$
where $H_{i j}=A_{i i} A_{j j}-A_{i j}^{2}$ and $K={\rm det} A$
is the Gauss-Kronecker curvature of $\mathcal{U}$.
Now, from \eqref{eqn:2.2} we have
$$2(H_{12}+H_{13}+H_{23})=\rho -(\kappa_1+\kappa_2)+(\kappa_1-\kappa_2)C.$$
By using  \eqref{eqn:3.1} and calculating the derivatives of the function ${\rm det} Q$ at $l=0$,
we  arrive at
$$
\begin{aligned}
\left(\frac{d\ {\rm det}\ Q}{d l}\right)(0)=&-H,\\
\left(\frac{d^{2} {\rm det}\ Q}{d l^{2}}\right)(0)
=&\rho-\frac{3(\kappa_1+\kappa_2)}{2}+\frac{\kappa_1-\kappa_2}{2}C,\\
\left(\frac{d^{4} {\rm det}\ Q}{d l^{4}}\right)(0)=&-4(1-C)\kappa_2
H_{12}-4(1+C)\kappa_1H_{13}-\frac{3C^2-2C-5}{4}\kappa_1^2-\frac{3C^2+2C-5}{4}\kappa_2^2\\
&+\frac{7+C^2}{2}\kappa_1\kappa_2-[(1+C)\kappa_1+(1-C)\kappa_2]\rho,\\
\left(\frac{d^{6} {\rm det}\ Q}{d l^{6}}\right)(0)=&
[ 6(1 - C)^2\kappa_2^2 +10(1-C^2)  \kappa_1\kappa_2]H_{12}+
[6 (1 + C)^2  \kappa_1^2  +10 (1 - C^2) \kappa_1\kappa_2]H_{13} \\
&+\frac{1}{4} [3(1+C)^2\kappa_1^2+10(1-C^2)\kappa_1\kappa_2+3(1-C)^2\kappa_2^2]\rho
+\frac{1}{8}\{(1 + C)^2 (5 C-7) \kappa_1^3 \\
&- (41 + 13 C - 17 C^2 +
11 C^3) \kappa_1^2\kappa_2 + (-41 + 13 C + 17 C^2 +
11 C^3) \kappa_1 \kappa_2^2 \\
&- (1 - C)^2 (7 +5 C) \kappa_2^3\}.
\end{aligned}
$$
Now we only discuss the following three cases for $\kappa_1\neq \kappa_2$ separately.

\noindent{\bf Case $\mathbb{S}^2\times \mathbb{H}^2$}. In this case, $\kappa_1=1$, $\kappa_2=-1$.
Taking into account that $\left(\frac{d^{2} {\rm det}\ Q}{d l^{2}}\right)(0)$,
$\left(\frac{d^{4} {\rm det}\ Q}{d l^{4}}\right)(0)$ and
$\left(\frac{d^{6} {\rm det}\ Q}{d l^{6}}\right)(0)$ are constant, we have
\begin{gather}
\label{eqn:3.2}
\rho+C=\alpha_1,\\
\label{eqn:3.3}
-1-2C^2+(4-4C)H_{12}-(4+4C)H_{13}-2C \rho=\alpha_2,\\
\label{eqn:3.4}
C+4C^3+(-4-12C+16 C^2)H_{12}
+(-4+12C+16 C^2)H_{13}+(4C^2-1)\rho =\alpha_3,
\end{gather}
where $\alpha_1, \alpha_2,\alpha_3$ are constant on $\mathcal{U}$.
				
Then, on $\mathcal{U}$, from \eqref{eqn:3.2}--\eqref{eqn:3.4} we obtain
\begin{gather}
\label{eqn:3.5}
\rho=\alpha_1-C,\\
\label{eqn:3.6}
H_{12}=-\frac{4 \alpha_1 C^2-(2 \alpha_1 -4\alpha_2-2 )C+\alpha_1-\alpha_2+\alpha_3-1}{8 (1-C)},\\
\label{eqn:3.7}
H_{13}=-\frac{4 \alpha_1 C^2+(2 \alpha_1 +4 \alpha_2 +2) C+\alpha_1+\alpha_2+\alpha_3+1}{8(1+C)}.
\end{gather}
				
In addition, in this case we need to derive $\left(\frac{d^{10} {\rm det}\ Q}{d l^{10}}\right)(0)$
to obtain another equation
$$
\begin{aligned}
\left(\frac{d^{10} {\rm det}\ Q}{d l^{10}}\right)(0)&=(128 C^4-80 C^3-96 C^2+40 C+8) H_{12}
+(128 C^4+80 C^3-96 C^2-40 C+8) H_{13}\\
&+(16 C^4  -12 C^2 +1)\rho+16 C^5-4 C^3-3 C.
\end{aligned}
$$
				
Then, substituting \eqref{eqn:3.5}--\eqref{eqn:3.7} into
$\left(\frac{d^{10} {\rm det}\ Q}{d l^{10}}\right)(0)=\alpha_4$ for a constant $\alpha_4$,
we get the following equation of $C$
$$
16 \alpha_2 C^3+(16 \alpha_1+12 \alpha_3) C^2+(4 \alpha_2+4) C-\alpha_1-2 \alpha_3-\alpha_4=0.
$$
Notice that the coefficients of $C$, $C^2$ and $C^3$ can not vanish simultaneously.
It follows that $C$ is constant on  $\mathcal{U}$.
Thus, by the continuity of the function $C$, we know that $C$ is constant on $\Sigma$.

\vskip 0.2cm
				
\noindent{\bf Case $\mathbb{S}^2\times \mathbb{R}^2$}. In this case, $\kappa_1=1$, $\kappa_2=0$.
From the expressions of $\left(\frac{d^2 {\rm det}\ Q}{d l^2}\right)(0)$,
$\left(\frac{d^{4} {\rm det}\ Q}{d l^{4}}\right)(0)$ and
$\left(\frac{d^{6} {\rm det}\ Q}{d l^{6}}\right)(0)$, we obtain
\begin{gather}
\label{eqn:3.8}
\frac{1}{2} (-3+C)+\rho=\alpha_1,\\
\label{eqn:3.9}
-\frac{1}{4} (1+C) \left(-5+3 C+4 \rho+16 H_{13}\right)=\alpha_2,\\
\label{eqn:3.10}
\frac{(1+C)^2}{8}(6\rho+48H_{13}+5C-7)=\alpha_3,
\end{gather}
where $\alpha_1$, $\alpha_2$ and $\alpha_3$ are constant on $\mathcal{U}$.
Then, from \eqref{eqn:3.8} and \eqref{eqn:3.9} we obtain
\begin{gather}
\label{eqn:3.11}
\rho=\alpha_1+\frac{3}{2}-\frac{C}{2},\\
\label{eqn:3.12}
H_{13}=-\frac{(1+ C)^2+4 \alpha_1(1+ C) +4 \alpha_2}{16 (1+C)}.
\end{gather}
Then, substituting \eqref{eqn:3.11} and \eqref{eqn:3.12} into \eqref{eqn:3.10}, we get
$$
(1+C)^3+6\alpha_1(1+C)^2+12\alpha_2(1+C)+8\alpha_3=0.
$$
It follows that $C$ is constant on $\mathcal{U}$.
Thus, by the continuity of the function $C$, we know that $C$ is constant on $\Sigma$.

\vskip 0.2cm
				
\noindent{\bf Case $\mathbb{H}^2\times \mathbb{R}^2$}. In this case,
$\kappa_1=-1$, $\kappa_2=0$. From the expressions of $\left(\frac{d^2 {\rm det}\ Q}{d l^2}\right)(0)$,
$\left(\frac{d^{4} {\rm det}\ Q}{d l^{4}}\right)(0)$ and
$\left(\frac{d^{6} {\rm det}\ Q}{d l^{6}}\right)(0)$,
we can get
\begin{gather}
\label{eqn:3.13}
\frac{1}{2} (3-C)+\rho=\alpha_1,\\
\label{eqn:3.14}
-\frac{1}{4} (1+C) \left(-5+3 C-4 \rho-16 H_{13}\right)=\alpha_2,\\
\label{eqn:3.15}
\frac{1}{8}(1+C)^2(6\rho+48H_{13}-5C+7)=\alpha_3,
\end{gather}
where $\alpha_1$, $\alpha_2$ and $\alpha_3$ are constant on $\mathcal{U}$.
Then, from \eqref{eqn:3.13} and \eqref{eqn:3.14} we obtain
\begin{gather}
\label{eqn:3.16}
\rho=\alpha_1-\frac{3}{2}+\frac{C}{2},\\
\label{eqn:3.17}
H_{13}=\frac{(1+ C)^2-4\alpha_1(1+ C) +4 \alpha_2}{16 (1+C)}.
\end{gather}
Then, substituting \eqref{eqn:3.16} and \eqref{eqn:3.17} into \eqref{eqn:3.15}, we get
$$
(1+C)^3-6\alpha_1(1+C)^2+12\alpha_2(1+C)-8\alpha_3=0.
$$
It follows that $C$ is constant on  $\mathcal{U}$.
Thus, by the continuity of the function $C$, we know that $C$ is constant on $\Sigma$.
Combined with the known results in \cite{Ur} and \cite{GMY},
we have completed the proof of Theorem \ref{thm:1.1}.
\qed

\vskip 0.2cm

\noindent{\bf Proof of Theorem \ref{thm:1.2}}.

Combining Theorem \ref{thm:1.1} and Theorem 1 of \cite{dos}, we know that
any isoparametric hypersurface $\Sigma$ in $M_{\kappa_1}^2\times M_{\kappa_2}^2$
for $\kappa_1\neq \kappa_2$ has constant principal curvatures.
Then, by applying Theorem 2 of \cite{dos}, we have completed the proof of
Theorem \ref{thm:1.2}.

\begin{acknow}
This is partially supported by National Natural Science Foundation of China
(Grant Nos. 11831005, 12061131014 and 12171437) and Research Ability Enhancement Program
for Young Teachers of Beijing University of Civil Engineering and Architecture (Grant No. X21025).
\end{acknow}
				

\end{document}